\documentclass[10pt]{article}

\newtheorem{theorem}{Theorem}

\newenvironment{proof}[1][Proof]{\noindent\textbf{#1.} }{\ \rule{0.5em}{0.5em}}

\begin{document}

\title{Forbidden Substrings on Weighted Alphabets}
\author{Amy N. Myers}
\maketitle

\begin{abstract}
In an influential 1981 paper, Guibas and Odlyzko constructed a generating function for the
number of length $n$ strings over a finite alphabet that avoid all
members of a given set of forbidden substrings. Here we
extend this result to the case in which the strings are weighted.
This investigation was inspired by the problem of counting
compositions of an integer $n$ that avoid all compositions of a
smaller integer $m$, a notion which arose from the
consideration of \emph{one-sided} random walks.
\end{abstract}

\section{Introduction}

In \cite{GO} Guibas and Odlyzko construct a generating function for the
number of length $n$ strings over a finite alphabet that avoid all
members of a given set of forbidden substrings. Here we
assign a weight to each letter of the alphabet, define the weight of a string to be
the sum of the weights of its letters, and determine a generating
function for the number of \emph{weight} $n$ strings that avoid a
particular set of forbidden substrings.  This investigation was
inspired by the problem of counting compositions of an integer $n$
(which can be viewed as weight $n$ strings over the alphabet
$\{1,2,3,\dots\}$) that do not contain a composition of a smaller
integer $m$ occurring in consecutive positions (\emph{i.e.}, avoid a
substring of weight $m$).  This latter problem arose from the
consideration of \emph{one-sided} random walks, which are introduced here, and
further investigated by Bender, Lawler, Pemantle, and Wilf in
\cite{BLPW}.  

Heubach and Kitaev have also extended Guibas and Odlyzko's results from words to compositions.  In \cite{HK} they consider length (number of parts) and weights in compositions over alphabets of the form $\{1,2,\dots,n\}$.  In this paper we consider arbitrary weighted alphabets.  For more on the combinatorics of compositions and words, see \cite{HM}.  

\section{Forbidden Substrings}

In this section we recall Guibas and Odlyzko's theorem concerning
forbidden substrings. 

A set $S=\{A,B,\dots ,T\}$ of strings over an alphabet $\Omega$ is
\textit{reduced} if no string contains another as a substring.  (In
particular, no string in $S$ is empty.) Let $f(n)$ denote the number
of length $n$ strings that avoid each member of $S$.  For each
string $H$ in $S$ let $f_{H}(n)$ denote the number of length $n$
strings that end with $H$ and avoid all members of $S$ except for
the single appearance of $H$ at the end.

Define generating functions $F(z)=\displaystyle\sum_{n\geq 0} f(n)/z^{n}$ and $F_{H}(z)=\displaystyle\sum_{n\geq 0}
f_{H}(n)/z^{n}.$

The \textit{correlation} of two strings $G$ and $H$, denoted $GH$,
is a string over $\{0,1\}$ with the same length as $G$.  The
$i^{th}$ character from the left in $GH$ is determined by placing
$H$ under $G$ so that the leftmost character of $H$ is under the
$i^{th}$ character (from the left) in $G$.  If all the pairs of
characters in the overlapping segment are identical, then the
$i^{th}$ character of $GH$ is 1.  If not, it is 0.  For example if
$\Omega =\{a,b\}$, $G=ababba$, and $H=abbab$, then $GH=001001$ as
illustrated below.

\begin{center}
$\begin{array}{ccccccccccc}
& a & b & a & b & b & a &  &  &  &  \\
0 & a & b & b & a & b &  &  &  &  &  \\
0 &  & a & b & b & a & b &  &  &  &  \\
1 &  &  & a & b & b & a & b &  &  &  \\
0 &  &  &  & a & b & b & a & b &  &  \\
0 &  &  &  &  & a & b & b & a & b &  \\
1 &  &  &  &  &  & a & b & b & a & b  \\
\end{array}$
\end{center}

Let $GH_{z}$ denote the correlation of $G$ and $H$ interpreted as a polynomial in the variable $z$.  In the
above example, $GH_{z}=z^{3}+1$.

\begin{theorem}
\label{GO} (Guibas, Odlyzko)  Given a reduced set $S=\{A,B,\dots
,T\}$ of strings over an alphabet of $q\geq 2$ characters, the
generating functions $F(z)$, $F_{A}(z)$, $F_{B}(z)$, $\dots$,
$F_{T}(z)$ satisfy the following system of linear equations:
\end{theorem}

$
\begin{array}{rrrrrrrrr}
(z-q)F(z) & + & zF_{A}(z) & + & \dots & + & zF_{T}(z) & = & z\\
F(z) & - & zAA_{z}F_{A}(z) & - & \dots & - & zTA_{z}F_{T}(z) & = & 0\\
F(z) & - & zAB_{z}F_{A}(z) & - & \dots & - & zTB_{z}F_{T}(z) & = & 0\\
\vdots & \vdots & \vdots & \vdots & \vdots & \vdots & \vdots & \vdots & \vdots\\
F(z) & - & zAT_{z}F_{A}(z) & - & \dots & - & zTT_{z}F_{T}(z) & = & 0\\
\end{array}
$

\bigskip

The fact that $S$ is reduced guarantees this system is nonsingular,
and we can solve for the generating functions as rational functions
of $z$.

\section{Weighted Strings}

Theorem \ref{GO} shows us how to construct a generating function for the number of
length $n$ strings that avoid each member of a given set of
forbidden substrings.  In this section we extend this result to
the case in which the strings are weighted, and count \emph{weight} $n$ strings.

A \textit{weighted} alphabet $w(\Omega)$ has each letter $h$
assigned a weight $w_h$. The \textit{weight} of a string
$H=h_{1}h_{2}\ldots h_{s}$ over $w(\Omega)$ is the sum
$w_H=\sum\limits_{i=1}^{t}w_{h_i}$ of the weights of the individual
letters.  A set $S$ of weighed strings is \textit{reduced} if no
string contains any other as a substring.

Given a set $S$ of reduced strings over a weighted alphabet
$w(\Omega)$, let $f(n)$ denote the number of weight $n$ strings that
do not contain any substring in $S$. Similarly for each $H$ in $S$
let $f_H(n)$ denote the number of weight $n$ strings that end with
$H$ and do not contain any substring in $S$ except for the single
appearance of $H$ at the end. Define $F(z)=\sum\limits_{n\geq 0}
f(n)/z^{n}$ and $F_{H}(z) = \sum\limits_{n\geq 0} f_{H}(n)/z^{n}$.
Note $f(0)=1$ counts the empty composition while $f_{H}(n)=0$ for
$n$ less than the number of letters in $H$.

Next we extend the notion of correlation for two strings to a
weighted version.  For the ordinary correlation $GH$ of two strings
$G$ and $H$, the $i^{th}$ character from the left is 1 if and only
if $G$ and $H$ overlap on the string $g_{i}g_{i+1}\ldots g_{r}$ for
some $r$.  The weighted correlation $w(GH)$ is a multiset, and the
weight $w_{g_{i}}+w_{g_{i+1}}+\ldots +w_{g_{r}}$ of the string on
which $G$ and $H$ overlap is in $w(GH)$.  More specifically, for any
two strings $G=g_{1}g_{2}\ldots g_{r}$ and $H=h_{1}h_{2}\ldots
h_{t}$ over a weighted alphabet $w(\Omega)$, the \textit{weighted
correlation} $w(GH)$ is a (possibly empty) mulitset.  This multiset
contains $k$ if and only if there is an $i$ such that
$h_{1}=g_{i},h_{2}=g_{i+1},\ldots$, $h_{r-i+1}=g_{r}$, and
$k=w_{h_{1}}+w_{h_{2}}+\ldots +w_{h_{r-i+1}}$ is the weight of the
overlap.

For example, let $w(\Omega)=\{1,2,\ldots \}$ with $w_i=i$. Set
$A=3,B=21,C=12,$ and $D=111$.  Then $w(AA)=w(BB)=w(CC)=\{3\},$
$w(CB)=\{2\},$ $w(BC)=w(DC)=w(BD)=\{1\},$ and $w(DD)=\{1,2,3\}.$ The
remaining weighted correlations are empty.  Note neither correlation
nor weighted correlation is commutative in general.

Finally we define $w(GH)_{z}$ to be the polynomial
$\sum\limits_{k\in w(GH)}z^{k}$.  When $w(GH)=\emptyset$, the
polynomial $w(GH)_{z}$ is 0.  Thus $w(DD)_{z}=z^{3}+z^{2}+z$, for
example.

We now prove an extension of Theorem \ref{GO}.

\begin{theorem}
\label{main} Given a reduced set $S=\{A,B,\ldots ,T\}$ of strings
over a weighted alphabet $w(\Omega)$, the generating functions
$F(z)$, $F_{A}(z)$, $F_{B}(z)$, $\ldots$, $F_{T}(z)$ satisfy the
following system of linear equations:
\end{theorem}

\noindent $
\begin{array}{ccccccc}
\frac{z-2}{z-1}F(z) & + F_{A}(z) & + F_{B}(z) & \cdots & + F_{T}(z) & = 1\\
F(z) & - w\left(AA\right)_{z}F_{A}(z) & -
w\left(BA\right)_{z}F_{A}(z) & \cdots &
- w\left(TA\right)_{z}F_{T}(z) & = 0\\
F(z) & - w\left(AB\right)_{z}F_{A}(z) & -
w\left(BB\right)_{z}F_{A}(z) & \cdots &
- w\left(TB\right)_{z}F_{T}(z) & = 0\\
\vdots & \vdots & \vdots & \vdots & \vdots & \vdots\\
F(z) & - w\left(AT\right)_{z}F_{A}(z) & -
w\left(BT\right)_{z}F_{A}(z) & \cdots &
-w\left(TT\right)_{z}F_{T}(z) & = 0\\
\end{array}
$

\bigskip

\begin{proof}  The first equation in the above system follows from the observation
that $f(n+1)+f_{A}(n+1)+\ldots+f_{T}(n+1)=f(n)+f(n-1)+\ldots+f(0)$.  This recurrence holds because any string $h_{1}h_{2}\ldots
h_{t}$ counted by one of $f(n+1)$, $f_{A}(n+1)$, $f_{B}(n+1)$,
$\ldots$, $f_{T}(n+1)$ arises by appending the character $h_{t}$ to
the string $h_{1}h_{2}\ldots h_{t-1}$ counted by $f(n+1-h_{t})$. The
right hand side of the recurrence equation is the coefficient of
$1/z^n$ in $z/(z-1)F(z)$, and the left hand side of the equation is
the coefficient of $1/z^n$ in
$z[F(z)-1]+zF_A(z)+zF_B(z)+\dots+zF_T(z)$. (Recall $f(0)=1$, but
$f_A(0)=f_B(0)=\dots f_T(0)=0$.)

The remaining equations result from the fact for any $H$ in $S$ we
have $f(n)=\sum\limits_{k\in w(AH)}f_{A}\left(n+k\right)
+\sum\limits_{k\in w(BH)}f_{B}\left(n+k\right) +\ldots
+\sum\limits_{k\in w(TH)}f_{T}\left(n+k\right)$.  To see this let $H=h_{1}h_{2}\ldots h_{t}$ and suppose
$Y=y_{1}y_{2}\ldots y_{s}$ is any string counted by $f(n).$  Let
$Z=z_{1}z_{2}\ldots z_{s+t}=y_{1}y_{2}\ldots y_{s}h_{1}h_{2}\ldots
h_{t}$ denote the concatenation of strings $Y$ and $H$.  Now $Z$
contains at least one string in $S$ as a substring.  Let
$G=g_{1}g_{2}\ldots g_{r}$ denote the leftmost such substring.  The
for some $u>s$ we have $g_{1}g_{2}\ldots g_{r}=z_{u-r+1}\ldots
z_{u-1}z_{u}$, and $z_{1}z_{2}\ldots z_{u}$ is counted by
$f_{G}(n+k)$ for some $k\in w\left(GH\right)$.

Conversely if $k\in w\left(GH\right)$, then any string counted by
$f_{G}\left(n\right)$ arises from the concatenation of a string $Y$
counted by $f(n)$ and $H$.  Thus the equality holds.  Since $\sum\limits_{n\geq 0}\sum\limits_{k\in w(GH)}
f_{G}\left(n+k\right)/z^{n}=\sum\limits_{k\in
w(GH)}z^{k}\sum\limits_{n\geq 0} f_{G}\left(n+k\right)
/z^{n+k}=w(GH)_{z}F_{G}(z)$, we obtain the remaining equations in the system.
\end{proof}

\bigskip

As was the case with Theorem \ref{GO}, the fact that $S$ is reduced
guarantees the system is nonsingular.  To see this, consider
the determinant\\

$\phi (z) = \det \left[
\begin{array}{ccccc}
\frac{z-2}{z-1} & 1 & 1 & \cdots & 1 \\
1 & -w\left(AA\right)_{z} & -w\left(BA\right)_{z} & \cdots & -w\left(TA\right)_{z} \\
1 & -w\left(AB\right)_{z} & -w\left(BB\right)_{z} & \cdots & -w\left(TB\right)_{z} \\
\vdots  & \vdots  & \vdots  & \ddots  & \vdots \\
1 & -w\left(AT\right)_{z} & -w\left(BT\right)_{z} & \cdots & -w\left(TT\right)_{z}\\
\end{array}
\right]$.

\ \\

\noindent Since $S$ is reduced the highest degree polynomial in each
column occurs on the diagonal.  When we expand $\phi (z)$, we have
$z-1$ in the denominator and a unique highest degree monomial
produced by the product of the diagonal terms in the numerator.  The
degree of this monomial is the sum $1+w_A+w_B+\dots +w_T$.  We can
therefore solve for $F(z),F_{A}(z),F_{B}(z),\dots ,F_{T}(z)$ and
find that each is a
rational function of $z$.\\

\section{Compositions}

The inspiration for extending Theorem \ref{GO} to Theorem \ref{main} came from the problem of
counting compositions of an integer $n$ that avoid compositions of a smaller integer $m$
occurring in consecutive positions.  For example, the composition
$2+4+1+1+4$ of $n=12$ contains the compositions $2+4$, $4+1+1$, and
$1+1+4$ of 6 in consecutive positions, while avoiding all
compositions of $m=3$ in consecutive positions.  (Note it does
contain the composition $2+1$ of 3 in nonconsecutive positions.)

We can apply Theorem \ref{main} to find, for example, a generating
function for the numbers of compositions of $n$ that avoid all
compositions of $m=3$ occurring in consecutive positions.  To do so
we view compositions as words over $w(\Omega) =\{1,2,3,\ldots \}$
with $w_i=i$. For example, we identify the composition $2+4+1+1+4$
as the word 24114. Set $S=\{A=3,B=21,C=12,D=111\}.$  The number of
compositions of $n$ that avoid all compositions of 3 occurring in
consecutive positions is given by the number of weight $n$ strings
over $w(\Omega)$ which do not contain any substring in $S$. Let
$f(n)$ denote this number. For each $H$ in $S$, let $f_{H}(n)$
denote the number of weight $n$ strings which end with $H$ and
contain no substring in $S$ except for the single occurrence of $H$
at the end.

Set $F(z)=\sum\limits_{n\geq 0} f(n)/z^{n}$ and
$F_{H}(z)=\sum\limits_{n\geq 0} f_{H}(n)/z^{n}$.  Earlier we recorded
the weighted correlation $w(GH)$ for each pair of strings in $S$.
We use this information to form the table below.  The polynomial
$w(GH)_{z}$ appears in row $H$ and column $G$.

\begin{center}
$\begin{array}{ccccc}
\ & A & B & C & D\\
A & z^{3} & 0 & 0 & 0 \\
B & 0 & z^{3} & z^2 & 0 \\
C & 0 & z & z^{3} & z \\
D & 0 & z & 0 & z^{3}+z^{2}+z\\
\end{array}$
\end{center}

Theorem \ref{main} guarantees the generating functions satisfy the
following system of equations:\\

$
\begin{array}{cccccc}
\frac{z-2}{z-1}F(z) & +F_{A}(z) & +F_{B}(z) &
+F_{C}(z) & +F_{D}(z) & =1\\
F(z) & -z^{3}F_{A}(z) & \ & \ & \ & =0\\
F(z) & \ & -z^{3}F_{B}(z) & -z^2F_{C}(z) & \ & =0\\
F(z) & \ & -zF_{B}(z) & -z^{3}F_{C}(z) & -zF_{D}(z) &  =0\\
F(z) & \ & -zF_{B}(z) & \ & -(z^{3}+z^{2}+z)F_{D}(z) & =0\\
\end{array}
$

\ \\

\noindent Solving this system yields\\

$
\begin{array}{lll}
F(z) & = &
(z^{8}-2z^{5}+z^{3})/(z^{8}-z^{7}-z^{6}+z^{5}-z^{4}-z^{3}-z^{2}+z+1)\\
\ & = &
1+1/z+2/z^2+2/{z^{4}}+3/{z^{5}}+9/{z^{6}}+12/{z^{7}}+20/{z^{8}}\ldots
\end{array}
$
 \ \\

\noindent as desired.

\section{Motivating Problem}

The question of composition avoidance arose from the consideration
of board games in which a roll of one or more (fair, 6-sided) dice determines the number of ``squares" a player moves forward on a given turn.  Some squares are undesirable to land on, and one
would like to know the probability of avoiding them, given the number of
squares separating a particular ``bad" square from one's current square.

To solve this problem, we replace the board with a finite number of squares, through which we cycle repeatedly, with an infinite succession of squares extending in one direction.  A sequence of dice rolls determines a \textit{one-sided random walk} which begins on square 0, and continues through squares 1, 2, 3, and so on, landing on some squares while avoiding others.  What is the probability that a one-sided random walk avoids square $m$?  

We record a one-sided random walk as an ``infinite composition'' of positive integer parts.  For example, $1+2+2+\dots$ indicates a sequence of rolls beginning with a roll of 1 followed by two rolls of 2.  What is the probability that a one-sided random walk avoids an initial composition of $m$?

Let $P(m)$ denote the probability that a one-sided random walk begins with a composition of $m$, and define $P(0)=1$.  In the simplest case, we use a single die to determine the size of each step in the walk, and compute $P(m)$ using the observation that $P(m)=\frac{1}{6}P(m-1)+\frac{1}{6}P(m-2)+\frac{1}{6}P(m-3)+\frac{1}{6}P(m-4)+\frac{1}{6}P(m-5)+\frac{1}{6}P(m-6)$.  From the recurrence we obtain the generating function
$$g(z) = \displaystyle\sum_{m \geq 0}P(m)z^m = \displaystyle\frac{1}{1-\frac{1}{6}{(z+z^2+z^3+z^4+z^5+z^6)}}$$
which converges for $|z|<1$.  Since $g(z)$ has a simple pole at $z=1$, and the residue there is $-\frac{2}{7}$, we know $P(m)\approx \frac{2}{7}$ for large $m$.  For large $m$, the probability that a one-sided random walk avoids square $m$ is therefore $1-P(m)\approx\frac{5}{7}$.  We can arrive at the same result using the fact that the recurrence for $P(m)$ has constant coefficients.  Specifically $P(m) = \frac{2}{7}+c_1r_1^m+c_2r_2^m+c_3r_3^m+c_4r_4^m+c_5r_5^m$, where $|r_i|<1$ for $1 \leq i \leq 5$.

\bigskip

More generally, we can consider one-sided random walks in which
$p_{i}$ is the probability of moving $i$ squares on a given turn.  When $p_{i}$ is determined by the roll of \textit{two} dice, we obtain $P(m)\approx \frac{1}{7}.$
The notion of one-sided random walks is considered further by Bender, Lawler, Pemantle, and Wilf in
\cite{BLPW}.    They compute, for instance, the probability of a ``collision" when two players take simultaneous one-sided random walks.  If $C(m)$ is the probability of a collision
for the first time on square $m$, then $$\sum\limits_{m \geq 0} C(m)x^{2m} = 1-\frac{1}{\frac{1}{2\pi}
\int\limits_{0}^{2\pi} \frac{d\theta}{1-|p(xe^{i\theta})|^{2}}}$$ where $p(z)=\sum\limits_{i\geq 1}p_{i}z^{i}.$

\bigskip

It is easy to count \emph{finite} compositions of an integer $n$ that avoid an initial composition of $m<n$.  Compositions that begin with an initial composition of $m$ have the form $\tau+\sigma$, where $\tau$ is a composition of $m$ and $\sigma$ is a composition of $n-m$.  There are $2^{m-1}\cdot2^{n-m-1} = 2^{n-2}$ compositions of $n$ that begin with a composition of $m$, and therefore also $2^{n-2}$ that avoid an initial composition of $m$, independent of our choice for $m$.  In other words, \emph{the
probability that a randomly selected composition of $n$ avoids an initial composition of $m < n$ is the same as the probability that it doesn't, for all such $m$}.

The mathematical literature contains numerous results concerning permutations and multiset permutations (which can be viewed as compositions) that avoid particular \emph{patterns}, \emph{i.e.}, permutations on fewer letters (see \cite{BonaPerms}, for example, for an introduction to the field).  The above investigation can be framed in this context as follows.  We know we can easily count compositions of $n$ that avoid an initial composition of $m<n$.  This results suggests the more general goal of  counting compositions of $n$ that avoid a composition of $m$ anywhere.  We can interpret this statement in several ways:
\begin{enumerate}
\item Count compositions of $n$ that avoid \textit{all} compositions of $m$ occurring in consecutive positions.
\item Count compositions of $n$ that avoid a \textit{particular} composition $\tau$ of $m$
occurring in consecutive positions.
\item Count compositions of $n$ that avoid a \textit{particular} composition $\tau$ of $m$ in
(possibly) nonconsecutive positions.
\item Count compositions of $n$ that avoid \textit{all} compositions of $m$ in (possibly) nonconsecutive positions.
\end{enumerate}

Here we have solved 1 and 2 with Theorem \ref{main}.  Problem 3 is straightforward, and 4 is open.

The problems above use the word ``avoid" in a narrow sense compared to that for patterns.  We can define a notion of composition avoidance analogous to that for pattern avoidance.  To do so we view the compositions of $n$ as multiset permutations.  For example, the compositions 1+1+2, 1+2+1,
and 2+1+1 correspond to the permutations 112, 121, and 211 of the
multiset $\{1^2,2\}$.  We identify the compositions of $n=4$ with
permutations of the multisets $\{4\}$, $\{1,3\}$, $\{2^2\}$,
$\{1^2,2\}$, and $\{1^4\}$.  It is well-known that the number of
permutations of a set of $n$ letters that avoid a pattern $\pi$ of 3
letters is independent of $\pi$.  Since the same result holds for
multisets (see \cite{MPA} or \cite{SW}), we see that the number of
compositions of $n$ that avoid a ``composition pattern" (I suggest
the term \emph{motif}) $\pi$ with 3 distinct parts is independent of
the parts. It would be interesting to investigate \emph{motif
avoidance} for other motifs.  The ($1+2$)-avoiding compositions are
the partitions.  How about the ($1+2+1$)-avoiding compositions?

\end{document}